\documentclass[11pt,a4paper]{article}

\usepackage{color}
\usepackage{graphics,epsfig}
\usepackage{amsfonts,amssymb,amsmath,amsthm}

\newtheorem{lem}{Lemma}[section]
\newtheorem{cor}[lem]{Corollary}
\newtheorem{thm}[lem]{Theorem}
\newtheorem{prop}[lem]{Proposition}
\theoremstyle{definition}
\newtheorem{defi}[lem]{Definition}
\theoremstyle{remark}
\newtheorem{rem}[lem]{Remark}

\numberwithin{equation}{section}

\renewcommand{\div}{\operatorname{div}}
\newcommand{\grad}{\operatorname{grad}}

\newcommand{\ep}{\varepsilon}
\newcommand{\ue}{u^\ep}

\newcommand{\n}{\nabla }
\newcommand{\om}{\Omega}
\newcommand{\emut}{e^{\mu t/\ep ^2}}

\newcommand{\ombar}{\overline{\Omega}}
\newcommand{\edeux}{\displaystyle{\frac{1}{\ep^2}}}
\newcommand{\R}{\mathbb{R}}

\newcommand{\vsp}{\vspace{8pt}}
\newcommand{\di}{\displaystyle}
\newcommand{\Pe}{(P^{\;\!\ep})}

\newcommand {\Q}{Q_T}
\newcommand{\regionunzero}{\om _0 ^{(1)}}
\newcommand{\regionzerozero}{\om _0 ^{(0)}}

\newcommand{\support}{\om _t ^{supp}}

\newcommand{\gammasupport}{\Gamma _t ^{supp}}
\newcommand{\gammasupportzero}{\Gamma _{t_0} ^{supp}}
\newcommand{\nusupport}{\nu _t ^{supp}}
\newcommand{\CG}{C^\star}

\title{Generation of interface for an Allen-Cahn equation \\
with nonlinear diffusion}
\author{ }
\date{}

\begin{document}

\maketitle \vspace{-20 mm}

\begin{center}

{\large\bf Matthieu Alfaro }\\[1ex]
I3M, Universit\'e de Montpellier 2,\\
CC051, Place Eug\`ene Bataillon, 34095 Montpellier Cedex 5, France,\\[2ex]

{\large\bf Danielle Hilhorst }\\[1ex]
CNRS et Laboratoire de Math\'ematiques\\
Universit\'e de Paris-Sud 11, 91405 Orsay Cedex, France. \\[2ex]

\end{center}

\vspace{15pt}


\begin{abstract}

In this note, we consider a nonlinear diffusion equation with a
bistable reaction term arising in population dynamics. Given a
rather general initial data, we investigate its behavior for small
times as the reaction coefficient tends to infinity: we prove a
generation
of interface property. \\

\end{abstract}

\section{Introduction}\label{intro-poreux}

We consider the degenerate parabolic problem
\[
 \Pe \quad\begin{cases}
 u_t=\Delta (u^m)+\edeux f(u)&\text{in }Q_T:=\om \times (0,T)\vspace{3pt}\\
 \di \frac{\partial (u^m)}{\partial \nu} = 0 &\text{on }\partial \om \times (0,T)\vspace{3pt}\\
 u(x,0)=u_0(x) &\text{in }\om\,,
 \end{cases}
\]
with $\ep >0$ a small parameter. Here  $\om$ is a smooth bounded
domain in $\R^N$ ($N\geq 2$), $\nu$ is the Euclidian unit normal
vector exterior to $\partial \om$ and $m\geq 2$.

We assume that $f$ is smooth and has exactly three zeros $0<a<1$
such that
\begin{equation}\label{der-f-poreux}
f'(0)<0\,, \qquad f'(a)>0\,, \qquad f'(1)<0\,.
\end{equation}

Moreover we suppose that the initial function $u_0 \in
C^2(\ombar)$ is nonnegative, with support $ Supp\,
u_0:=\{x\in\om|\ u_0(x)>0\} \subset\subset \om$. Furthermore we
define the initial interface $\Gamma _0$ by
$$
\Gamma _0:=\{x\in\om|\ u_0(x)=a \}\,,
$$
and suppose that $\Gamma _0$ is a smooth hypersurface without
boundary such that
\begin{equation}\label{dalltint-poreux}
\Gamma _0 \subset\subset \Omega \quad \mbox { and } \quad \n
u_0(x) \neq 0\quad\text{if $x\in\Gamma _0$\,,}
\end{equation}
\begin{equation}\label{initial-data-poreux}
u_0>a \quad \text { in } \quad \regionunzero\,,\quad u_0<a \quad
\text { in } \quad \regionzerozero\,,
\end{equation}
where $\regionunzero$ denotes the region enclosed by $\Gamma _0$
and $\regionzerozero$ the region enclosed between $\partial \om$
and $\Gamma _0$.\\

We prove a generation of interface property, namely that the
solution $\ue$ quickly becomes close to $1$ or $0$, except in an
$\mathcal O(\ep)$ neighborhood of the initial interface $\Gamma
_0$, creating a steep transition layer around $\Gamma _0$. More
precisely, we are in presence of the following phenomenon: in the
very early stage, the nonlinear diffusion term is negligible when
compared with the reaction term $\ep ^{-2}f(u)$. Hence, under the
rescaling in time $\tau=t/\ep^2$, the equation is well
approximated by the ordinary differential equation $u_\tau=f(u)$.
In view of the bistable nature of $f$, $\ue$ quickly approaches
the stable equilibria of the ordinary differential equation, $0$
or $1$, and an interface is formed between the regions
$\{\ue\approx 0\}$ and
$\{\ue\approx 1\}$.\\

The organization of this note is as follows. In Section
\ref{s:biology} we briefly explain how Problem $(P^\varepsilon)$
arises in population dynamics. In Section \ref{s:existence}, we
recall known results about the well-posedness of Problem $\Pe$ and
a comparison principle. In Section \ref{s:generation}, we prove
the generation of interface property for Problem $\Pe$. To that
purpose we construct sub- and super-solutions by modifying the
solution of the corresponding ordinary differential equation
$u_t=\ep ^{-2}f(u)$. We also show the optimality of the generation
time ${t}^{\,\ep}:=f'(a) ^{-1}\ep^2|\ln\ep|$ and prove that the
thickness of the interface is of order $\mathcal O(\ep)$ at the
generation time ${t}^{\,\ep}$. Our method of proof follows the
same lines as that of \cite{AHM} and \cite{A}. It is slightly
different from those of Xinfu Chen \cite{C1} and \cite{C2}, who
transforms the reaction function $f$. We postpone to future work
the study of the interface motion after the generation
time of the interface. \\

Finally let us mention articles involving the singular limit of
reaction-diffusion equations with nonlinear diffusion. Feireisl
\cite{F} studies the singular limit of a degenerate parabolic
equation in the whole space $\R ^N$. He studies the problem in the
scaling
\begin{equation}\label{feireisl}
u_t=\ep \Delta (u^m)+\frac 1 \ep f(u)\,,
\end{equation}
where, in the limit $\ep \to 0$, the limit free boundary moves
according to motion by constant speed (that of a related traveling
wave). In a similar scaling, Hilhorst, Kersner, Logak and Mimura
\cite{HKLM} investigate the singular limit of this equation in a
bounded domain with a monostable reaction term. In both of these
papers, they prove that the solution $\ue$ of the nonlinear
diffusion equation converges to 0 or 1 on both sides of an
interface moving with constant normal velocity. In this scaling,
proofs about the interface motion can be performed with using only
one term in the asymptotic expansion whereas we would need to use
two terms in the case of Problem $(P^\varepsilon)$ as well as a
suitable linearization procedure; this is far from trivial here
since Problem $\Pe$ is degenerate parabolic.

\section{The biological context}\label{s:biology}
In this section, we discuss nonlinear diffusion in population
dynamics. It is well-known that the control of a population can be
achieved by introducing density dependent birth or death rates. In
\cite{GN}, Gurney and Nisbet show that the introduction of a
nonlinearity into the dispersal behavior of a species --- which
behaves in an otherwise linear way--- can, in an inhomogeneous
environment, lead to a regulatory effect. More precisely, they
consider the equation
$$
u_t=-\div {\bf j}+G(x)u\,,
$$
where $u(x,t)$ denotes the population density, $G=G(x)$ the growth
function only depending on the location and ${\bf j}(x,t)$ the
local population current density. By using the well-known random
motion model one obtains the linear equation $u_t=\Delta u+G(x)u$.
Another possibility is to choose the biased random motion model
where movements are largely random but slightly modified by the
distribution of the fellows; the corresponding equation is then
written as $u_t=\Delta u+\div(u\grad u)+G(x)u$. Nevertheless, Carl
\cite{Carl} has observed that arctic ground squirrels migrate from
densely populated areas into sparsely populated ones, even when
the latter is less favorable (burrow sites not available,
intensive predation). For such species, migration to avoid
crowding, rather than random motion, is the primary cause of
dispersal. To describe such movements, Gurney and Nisbet use the
directed motion model where individuals can only stay put or move
down the population gradient; this model yields the degenerate
parabolic equation
\begin{equation}\label{gurney}
u_t=\Delta (u^2)+G(x)u\,.
\end{equation}
In \cite{GN}, the authors perform a qualitative analysis of the
three different dispersal models (random motion, biased random
motion and directed motion). They conclude that the introduction
of density dependent diffusion is efficient to study the dynamics
of a population which regulates its size below the carrying
capacity set by the supply of nutrients. \\

Gurtin and Mac Camy \cite{GM} proposed the class of equations
which we study here and which involves degenerate diffusion and
nonlinear reaction, namely
\begin{equation}\label{general}
u_t=\Delta (u^m)+f(u)\,,\quad \quad m\geq 2\,.
\end{equation}
In absence of a reaction term, equation \eqref{general} reduces to
the so-called porous medium equation
\begin{equation}\label{pme}
u_t=\Delta (u^m)\,,
\end{equation}
which describes, among others, the flow of an ideal gas in a
homogeneous medium ($m\geq 2$), groundwater infiltration
(Boussinesq's equation, $m=2$), the spread of a thin viscous film
under gravity ($m=4$), and thermal propagation in plasma ($m
\simeq 6$). The porous medium equation has been extensively
investigated in literature: we refer to the book of V\'asquez
\cite{V} and the references therein. The main feature of these
equations is that they degenerate at the points where $u=0$. As a
consequence, a loss of regularity of solutions occurs and
disturbances propagate with finite speed. This phenomenon
contrasts with the infinite speed of propagation of solutions of
the
heat equation $u_t=\Delta u$.  \\

\section{Comparison principle and well-posedness}\label{s:existence}
Since the diffusion term degenerates at the points where
$u^\varepsilon = 0$, $u^\varepsilon$ is not smooth. This leads us
to define a notion of weak solution for Problem $\Pe$, in a
similar way as it is done by Aronson, Crandall and Peletier
\cite{ACP} for a corresponding one-dimensional problem.
\begin{defi}\label{definition-weaksol-poreux}
A function $u^\varepsilon : [0,\infty)\to L^1(\om)$ is a weak
solution of Problem $(P^\varepsilon)$ if, for all $T>0$,
\begin{enumerate}
\item $u^\varepsilon \in C\left([0,\infty);L^1(\om)\right)\cap
L^\infty (\Q)$\,; \item $u^\varepsilon$ satisfies the integral
equality
\begin{multline}\label{deqexi-poreux}
\int_ \om
{u^\varepsilon}(T)\varphi(T)-\int\int_{\Q}({u^\varepsilon}\varphi
_t +(u^\varepsilon)^m\Delta \varphi)= \int _\om u_0\varphi(0)
\\+\int\int _{\Q} \edeux f(u^\varepsilon)\varphi\,,
\end{multline}
for all $\varphi \in C^2(\overline{\Q})$ such that $\varphi \geq
0$ and $\di \frac {\partial \varphi}{\partial \nu}=0$ on $\partial
  \om$.
\end{enumerate}
A sub-solution (respectively a super-solution) of Problem
$(P^\varepsilon)$ is a function satisfying (i) and (ii) with
equality replaced by $\leq$ (respectively $\geq$).
\end{defi}

\begin{thm}[Existence and
comparison principle]\label{Existence-comparison} Let $T>0$ be
arbitrary. The following properties hold:
\begin{enumerate}
 \item Let $u^-$ and $u^+$ be a sub-solution and a super-solution
of Problem $(P^\varepsilon)$ with initial data $u_0^-$ and $u_0^+$
respectively.
$$
\mbox{If~}\quad u_0^- \leq u_0^+\,\quad \mbox{~then~}\quad u^-\leq
u^+ \mbox{~in~} Q_T\,;
$$
 \item Problem $(P^\varepsilon)$ has a unique weak solution $u^\varepsilon$
which is such that
\begin{equation}\label{encadrement}
0\leq u^\varepsilon \leq \max(1,\Vert u_0 \Vert _{L^\infty(\om)})
\text{~in~} Q_T\,;
\end{equation}
\item $u^\varepsilon \in C(\overline{Q_T}).$
\end{enumerate}
\end{thm}

The proof of Theorem \ref{Existence-comparison} is standard; it
can be performed by using the same lines as that of Theorem 5 in
\cite{ACP}. The continuity of $u^\varepsilon$ follows from
\cite{DB}.\\

The following result turns out to be an essential tool when
constructing smooth sub- and super-solutions of Problem $\Pe$.

\begin{lem}\label{lemma-sup-poreux}
Let $u^\varepsilon$ be a continuous nonnegative function in
$\ombar \times [0,T]$. Define $\support = \{x \in \om |\
u^\varepsilon (x,t)>0\}$ and $\gammasupport=\partial \support$ for
all $t\in[0,T]$. Suppose that the family $\Gamma:=\cup _{0< t \leq
T} \gammasupport  \times \{t\}$ is sufficiently smooth and let
$\nusupport$ be the outward normal vector on $\gammasupport$.
Suppose moreover that
\begin{enumerate}
\item $\nabla (u^\varepsilon)^m \mbox{~is~continuous~in~}
\overline{\Omega} \times [0,T]\,;$ \item ${\cal
L}[u^\varepsilon]:=u^\varepsilon_t -\Delta ({u^\varepsilon})^m
-\edeux f(u^\varepsilon)
= 0 \mbox{ in } \{(x,t) \in \ombar \times [0,T] \\
\mbox{~such~that~} {u^\varepsilon}(x,t) > 0\}\,;$ \item
$\di{\frac{\partial (u^\varepsilon)^m}{\partial \nusupport}}=0\;
\text{ on }\
\partial \support,\; \text{ for all }\ t\in[0,T]$\,.
\end{enumerate}
Then $u^\varepsilon$ is a solution of Problem $(P^\varepsilon)$.
Similarly $u$ is a sub-solution (respectively a super-solution) of
Problem $(P^\varepsilon)$ if the equality in (ii) is replaced by
$\leq$ (resp. $\geq$) and if the equality in (iii) is replaced by
$\leq$ (resp. $\geq$).
\end{lem}
We refer to \cite{HKLM} for the proof.

\section{Generation of interface}
\label{s:generation}

In this section we prove that, given a nearly arbitrary initial
function $u_0$, the solution $\ue$ quickly becomes close to $1$ or
$0$, except in an $\mathcal O(\ep)$ neighborhood of the initial
interface $\Gamma_0$, creating a steep transition layer around
$\Gamma _0$. The time needed to develop such a transition layer is
of order $\mathcal O (\ep^2|\ln\ep|)$.

\begin{thm}[Generation of interface]\label{g-th-gen-poreux}
Assume $m\geq 2$. Let $\gamma \in (0,\min (a,\\
1 -a))$ be arbitrary and define $\mu$ as the derivative of $f(u)$
at the unstable equilibrium $u=a$, namely
\begin{equation}\label{g-def-mu-poreux}
\mu=f'(a)\,.
\end{equation}
Moreover, set
$$
t^{\,\ep}:=\mu ^{-1}\ep ^2|\ln\ep|\,.
$$
Then there exist positive constants $\ep_0$ and $M_0$ such that,
for all $\ep\in(0,\ep _0)$,
\begin{enumerate}
\item for all $x \in \om$, we have that
$$
0 \leq u^\ep(x,t^{\,\ep}) \leq 1+\gamma\,;
$$
\item for all $x\in\om$ such that $|u_0(x)-a|\geq M_0 \ep$, we
have that
\begin{align}
&\text{if}\;~~u_0(x)\geq
a+M_0\ep\;~~\text{then}\;~~u^\ep(x,t^{\,\ep})
\geq 1-\gamma\,,\label{g-part2}\\
&\text{if}\;~~u_0(x)\leq
a-M_0\ep\;~~\text{then}\;~~u^\ep(x,t^{\,\ep}) \leq \gamma
\label{g-part3}\,.
\end{align}
\end{enumerate}
\end{thm}

Theorem \ref{g-th-gen-poreux} will be proved by constructing a
suitable pair of sub- and super-solutions. As mentioned above, the
nonlinear diffusion term is negligible in this early stage so that
the behavior of the solution is governed by the ordinary
differential equation $u_t=\edeux f(u)$.

An immediate consequence of Theorem \ref{g-th-gen-poreux} is that
$\ue(x,t^{\,\ep})$ is close to 0 or 1, except in an $\mathcal
O(\ep)$ neighborhood of the initial interface $\Gamma _0$. In
other words, the transition layers which have developed have an
$\mathcal O(\ep)$ thickness.

\begin{cor}[Thickness of the transition layers at time $t^{\,\ep}$]\label{cor-thi}
Let $\eta \in (0,\min (a, 1 -a))$ be an arbitrary constant. Then
there exist positive constants $\ep _0 $ and $\mathcal C$ such
that, for all $\ep \in(0,\ep _0)$,
\begin{equation}\label{resultat}
\ue(x,t^{\,\ep}) \in
\begin{cases}
\,[0,1+\eta]&\quad\text{if}\quad
x\in\mathcal N_{\mathcal C\ep}(\Gamma_0)\\
\,[0,\eta]&\quad\text{if}\quad x\in \om ^{(0)} _0
\setminus\mathcal N_{\mathcal C\ep}(\Gamma
_0)\\
\,[1-\eta,1+\eta]&\quad\text{if}\quad x\in\om ^{(1)}
_0\setminus\mathcal N_{\mathcal C\ep}(\Gamma _0)\,,
\end{cases}
\end{equation}
where $\om ^{(1)} _0$ denotes the region enclosed by $\Gamma _0$,
$\om ^{(0)} _0$ the region enclosed between $\partial \om$ and
$\Gamma_0$, and
$$
\mathcal N _r(\Gamma _0):=\{x\in \om,\, dist(x,\Gamma _0)<r\}
$$
denotes the $r$-neighborhood of $\Gamma _0$.
\end{cor}

We will also show that the generation time
${t}^{\,\ep}:=\mu^{-1}\ep^2|\ln\ep|$ is optimal. In other words,
the interface is not fully developed until $t$ becomes close to
$t^{\,\ep}$. More precisely, the following result holds.

\begin{prop}\label{pr:optimal-time}
Denote by $t^{\,\ep}_{min}$ the smallest time such that
\eqref{resultat} holds. Then there exists a constant $b=b(\mathcal
C)$ such that
\[
{t}^{\,\ep}_{min}\geq \mu^{-1} \ep^2 (|\ln\ep| - b)\,
\]
for all $\,\ep \in (0,\ep_0).$
\end{prop}

\subsection{Proof of the generation of interface property}

\subsubsection{The bistable ordinary differential equation}

Let us first consider the problem without diffusion, namely
\begin{equation*}\label{no-diffusion}
\bar{u}_t=\frac{1}{\ep^2}\,f(\bar{u})\,, \qquad
\bar{u}(x,0)=u_0(x)\,.
\end{equation*}
Its solution can be written in the form
\[
\bar{u}(x,t)=Y\left(\frac{t}{\ep^2},\,u_0(x)\right)\,,
\]
where $Y(\tau,\xi)$ is the solution of the ordinary differential
equation
\begin{equation}\label{ode-poreux}
\left\{\begin{array}{ll} Y_\tau (\tau,\xi)&=f(Y(\tau,\xi)) \quad
\text { for } \tau >0 \vspace{3pt}\\
Y(0,\xi)&=\xi\,.
\end{array}\right.
\end{equation}
Here $\xi$ ranges over the interval $(-C_0,C_0)$, where
$C_0:=\Vert u_0 \Vert _{L^\infty(\om)} +1$. We claim that $Y$ has
the following properties.

\begin{lem}\label{properties-Y}There exists a positive constant $C$ such that the following
holds
\begin{enumerate}
\item
$\text{If }\ \xi >0 \;\text{ then }\ Y(\tau,\xi)>0\,,$\\
$\text{If }\ \xi <0 \;\text{ then }\ Y(\tau,\xi)<0\,;$ \item
$|Y(\tau,\xi)|\leq C_0\,;$ \item $Y _ \xi (\tau,\xi) > 0\,;$ \item
$|\di{\frac{Y_{\xi\xi}}{Y_\xi}(\tau,\xi)}|\leq C (e^{\mu
\tau}-1)$\,,
\end{enumerate}
for all $\tau >0$ and all $\xi \in(-C_0,C_0)$.
\end{lem}
Properties (i) and (ii) are direct consequences of the profile of
$f$ --- more precisely of the sign conditions $f>0$ in
$(-\infty,0)\cup(a,1)$ and $f<0$ in $(0,a)\cup (1,\infty)$--- and
of the qualitative properties of the solution of the bistable
ordinary differential equation \eqref{ode-poreux}; for proofs of
(iii) and (iv) we refer to \cite{AHM}, subsection 3.1.\qed

\subsubsection{Construction of sub- and super-solutions}

We use the notation $a^+=\max(a,0)$. The sub- and super-solutions
are given by
\begin{equation}\label{w+-}
w_\ep^\pm(x,t)=\left[Y\left(\frac{t}{\ep^2},\,u_0(x)\pm\ep^2\CG(\emut-1)\right)\right]^+\,.
\end{equation}

\begin{lem}\label{g-w}
There exist positive constants $\ep_0$ and $\CG$ such that, for
all $\, \ep \in (0,\ep _0)$, $(w_\ep^-,w_\ep^+)$ is a pair of sub-
and super-solutions for Problem $\Pe$, in the domain $\ombar\times
[0,\mu ^{-1} \ep^2|\ln \ep|]$. Moreover, since also
$w^-_\ep(x,0)=w^+ _\ep(x,0)=u_0(x)$, it follows that
\begin{equation}\label{g-coincee1}
w_\ep^-(x,t) \leq u^\ep(x,t) \leq w_\ep^+(x,t) \quad\
\mbox{~for~all~} (x,t) \in \ombar\times [0,\mu ^{-1} \ep^2|\ln
\ep|]\,.
\end{equation}
\end{lem}

{\noindent \bf Proof.} In order to prove that $(w_\ep ^-,w_\ep
^+)$ is a pair of sub- and super-solutions for Problem $\Pe$ for a
suitable choice of $\ep _0$ and $C_g$, we check that the
sufficient conditions in Lemma \ref{lemma-sup-poreux} are
satisfied.

As for the sub-solution $w_\ep ^-$, we remak that property (i) in
Lemma \ref{properties-Y} implies that, for all $t>0$,
$$
\begin{array}{ll}
\support [w_\ep^-]=\{x\in\om| \ u_0(x)>\ep ^2 C_g (\emut
-1)\}\vspace{3pt}\\
\gammasupport [w_\ep^-]:=\partial \support [w_\ep ^-]= \{x\in\om|
\ u_0(x)=\ep ^2 C_g (\emut -1)\}\,.
\end{array}
$$
Choose $(x_0,t_0)$ such that $x_0 \in\gammasupportzero [w_\ep
^-]$; for $(x,t)$ such that $x\in \support [w_\ep ^-]$ we have
$$
\n (w_\ep ^-)^m (x,t)=
mY^{m-1}Y_\xi\left(\frac{t}{\ep^2},\,u_0(x)-\ep^2\CG(\emut-1)\right)\n
u_0(x)\,.
$$
Since $Y(\tau,0)=0$ the equality above implies
$$
\lim _{\substack{(x,t)\to(x_0,t_0)\vspace{3pt}\\ x\in \support
[w_\ep ^-]}} \n (w_\ep ^-)^m (x,t)=0.
$$
Therefore conditions (i) and (iii) of Lemma \ref{lemma-sup-poreux}
are satisfied by the sub-solution.

As for the super-solution $w_\ep ^+$, we remark that property (i)
of Lemma \ref{properties-Y} implies that, for all $t>0$,
$$
\begin{array}{ll}
\support[w_\ep^+]=\om\vspace{3pt}\\
\gammasupport[w_\ep ^+]:=\partial \support [w_\ep ^+]=\partial
\om\,.
\end{array}
$$
Hence condition (iii) of Lemma \ref{lemma-sup-poreux} for the
super-solution is a direct consequence of the fact that $ Supp\,
u_0 \subset\subset \om$, whereas condition (i) is obviously
satisfied.

It remains to prove that
$$
{\cal L} [w_\ep ^-]:=(w_\ep ^-)_t-\Delta ({w_\ep ^-})^m-\edeux
f(w_\ep ^-)\leq 0\,,
$$
in $\{(x,t)\in \ombar\times [0,\mu ^{-1}\ep ^2 |\ln \ep|]
\mbox{~such~that~} w_\ep ^-(x,t)>0\}$ and that ${\cal L} [w_\ep
^+] \geq 0$ in $\{(x,t)\in \ombar\times [0,\mu ^{-1}\ep ^2 |\ln
\ep|]\}$.

In view of the ordinary differential equation \eqref{ode-poreux},
straightforward computations yield
\begin{multline}
{\cal L} [{w_\ep^-}] =-Y_\xi\Big[\CG \,\mu\, \emut+m(m-1)
Y^{m-2}Y_\xi
|\n u_0|^2\\
+mY^{m-1}{\frac{Y_{\xi\xi}}{Y_\xi}}\,|\nabla u_0|^2
+mY^{m-1}\Delta u_0\Big]\,,
\end{multline}
in $\support [w_\ep^-]$, where the function $Y$ and its
derivatives are taken at the point $(\tau,\xi)=(t /{\ep ^2},
u_0(x)-\ep^2\CG(\emut-1))$. Moreover since the term $m(m-1)
Y^{m-2}Y_\xi^2 |\n u_0|^2$ is nonnegative, it follows that
$$
{\cal L} [w_\ep^-] \leq -Y_\xi\Big[\CG \,\mu\,
\emut+mY^{m-1}{\frac{Y_{\xi\xi}}{Y_\xi}}\,|\nabla u_0|^2
+mY^{m-1}\Delta u_0\Big]\,.
$$
We note that, in the range $0 \leq t \leq \mu ^{-1} \ep ^2|\ln
\ep|$, we have, for $\ep _0$ sufficiently small,
$$
\xi=u_0(x)- \ep ^2 \CG(\emut-1)\in (-C_0,C_0)\,.
$$
We deduce from the properties (ii)-(iv) stated in Lemma
\ref{properties-Y} that there exist positive constants $C_1$ and
$C_2$ --- only depending on $m$, $C_0$, $C$, $\Vert \n u_0\Vert
_{L^\infty(\om)}$ and $\Vert \Delta u_0 \Vert _{L^\infty(\om)}$---
such that
$$
{\cal L}[w_\ep^-] \leq -Y_\xi \Big[(\CG\,\mu-C_1)\emut-C_2\Big]\,,
$$
which implies that ${\cal L} [w_\ep ^-] \leq 0$ if $\CG$ is chosen
large enough.

As for the super-solution we obtain
\begin{multline}
{\cal L} [w_\ep^+] =Y_\xi\Big[\CG \,\mu\, \emut-m(m-1)
Y^{m-2}Y_\xi
|\n u_0|^2\\
-mY^{m-1}{\frac{Y_{\xi\xi}}{Y_\xi}}\,|\nabla u_0|^2
-mY^{m-1}\Delta u_0\Big]\,,
\end{multline}
and the assumption that $m\geq 2$ gives an upper bound for
$|Y^{m-2}|$. Following the same argument as above one can prove
that ${\cal L} [w_\ep ^+] \geq 0$ for $\CG$ sufficiently large.
This completes the proof of Lemma \ref{g-w}.\qed

\subsubsection{Proof of Theorem \ref{g-th-gen-poreux}}

In order to prove Theorem \ref{g-th-gen-poreux} we first present
basic estimates of the function $Y$ after a time of order
$\tau\sim |\ln \ep|$.

\begin{lem}\label{after-time}
Let $\gamma \in (0,\min (a,1 -a))$ be arbitrary. There exist
positive constants $\ep_0$ and $C_Y$ such that, for all
$\,\ep\in(0,\ep _0)$,
\begin{enumerate}
\item for all $\xi\in (-C_0,C_0)$,
\begin{equation}\label{g-part11}
-\gamma \leq Y(\mu ^{-1} | \ln \ep |,\xi) \leq 1+\gamma\,;
\end{equation}
\item for all $\xi\in (-C_0,C_0)$ such that $|\xi-a|\geq C_Y \ep$,
we have that
\begin{align}
&\text{if}\;~~\xi\geq a+C_Y \ep\;~~\text{then}\;~~Y(\mu ^{-1}| \ln
\ep |,\xi)
\geq 1-\gamma\label{g-part22}\\
&\text{if}\;~~\xi\leq a-C_Y \ep\;~~\text{then}\;~~Y(\mu ^{-1}| \ln
\ep |,\xi)\leq \gamma \label{g-part33}\,.
\end{align}
\end{enumerate}
\end{lem}

These estimates illustrate the stability of the equilibria 0 and 1
for the bistable ordinary differential equation
\eqref{ode-poreux}. For more details we refer the reader to the
proof of Lemma 3.9 in \cite{AHM}.\qed

\vskip 8pt We are now ready to prove Theorem
\ref{g-th-gen-poreux}. By setting $t=\mu ^{-1} \ep ^2|\ln \ep|$ in
\eqref{g-coincee1}, we obtain
\begin{multline}\label{g-gr}
\left[Y\left(\mu ^{-1}|\ln \ep|, u_0(x)-(\CG \ep -\CG \ep ^2 )\right)\right] ^+\\
\leq u^\ep(x,\mu ^{-1} \ep^2|\ln \ep|) \leq \left[Y\left(\mu
^{-1}|\ln \ep|, u_0(x)+\CG \ep -\CG \ep ^2\right)\right]^+\,.
\end{multline}
Since, for $\ep_0$ small enough, $u_0(x)+ (\CG \ep -\CG \ep ^2)
\in (0,C_0)$, the assertion $(i)$ of Theorem \ref{g-th-gen-poreux}
is a direct consequence of \eqref{g-part11} and \eqref{g-gr}.

Next we prove \eqref{g-part2}. We choose $M_0$ large enough so
that $M_0\ep- \CG \ep+\CG \ep ^2 \geq C_Y \ep$.  Then, for all
$x\in \om$ such that $u_0(x)\geq a+M_0 \ep$, we have that $
u_0(x)-(\CG \ep -\CG \ep ^2) \geq a+C_Y \ep$, which we combine
with \eqref{g-gr} and \eqref{g-part22} to deduce that
\[
u^\ep(x,\mu^{-1} \ep ^2 | \ln \ep |)\geq 1-\gamma.
\]
The inequality \eqref{g-part3} can be shown in a similar way.
This completes the proof of Theorem \ref{g-th-gen-poreux}.\qed

\begin{rem}\label{gpaszero} Theorem \ref{g-th-gen-poreux}
remains true if we perturb the reaction function $f(u)$ by order
$\ep$, setting for instance ${\widetilde{f}}(u) = f(u)-\ep
g(x,t,u)$. To deal with this more general case, we proceed as
follows. We first consider a slightly perturbed reaction function,
namely ${f_\delta}(u)=f(u)+\delta$ which, for $\delta$ small
enough, is still of bistable type. Define $a(\delta)$ as its
unstable zero and $\mu(\delta):=f'(a(\delta))$ as the slope of
$f_\delta$ in this point. We then define $Y(\tau,\xi;\delta)$ as
the solution of the initial value problem
\begin{equation}
\left\{\begin{array}{ll} Y_\tau (\tau,\xi;\delta)&=f_\delta
(Y(\tau,\xi;\delta)) \quad \text { for }\tau >0\vspace{3pt}\\
Y(0,\xi;\delta)&=\xi\,.
\end{array}\right.
\end{equation}
Finally we construct a pair of sub- and super-solutions in the
form
\[
w_\ep^\pm(x,t)=\left[Y\left(\frac{t}{\ep^2},u_0(x)\pm\ep^2r(\pm
\ep \mathcal G,\frac{t}{\ep^2});\pm \ep \mathcal
G\right)\right]^+\,,
\]
where $r(\delta,\tau):=\CG(e^{\mu (\delta) \tau}-1)$ and $\mathcal
G:=\Vert g\Vert _{L^\infty(0,C_0)}$. For more details and proofs
we refer to \cite{AHM}, more precisely to Section 4 which deals
with a generation of interface property for an equation with
linear diffusion and an unbalanced reaction term $f(u)$.\qed
\end{rem}

\subsection{Proof of the optimality of the generation time}

We show below that the generation time
${t}^{\,\ep}:=\mu^{-1}\ep^2|\ln\ep|$ is optimal. In other words,
the interface is not fully developed before $t$ is close to
$t^{\,\ep}$. We will need the following lemma about the solution
$Y(\tau,\xi)$ of the corresponding bistable ordinary differential
equation.

\begin{lem}\label{lemme-opt}
Let $\eta \in (0,\min(a,1-a))$ be arbitrary. Then there exist
positive constants $C_1=C_1(\eta)$ and $C_2=C_2(\eta)$ such that,
if $\xi\in (a,1-\eta)$ then, for every $\tau>0$ such that
$Y(\tau,\xi)$ remains in the interval $(a,1-\eta)$, we have
\begin{equation}\label{g-est-Y-1}
C_1e^{\mu \tau}(\xi-a)\leq Y(\tau,\xi)-a \leq C_2e^{\mu
\tau}(\xi-a)\,.
\end{equation}
\end{lem}

We refer to \cite{AHM} Corollary 3.5 for the proof of Lemma \ref{lemme-opt}.\qed \\

{\noindent \bf Proof of Proposition \ref{pr:optimal-time}.} For
each $b>0$, we set
\[
t{\, ^\ep} (b):=\mu ^{-1} \ep ^2 (|\ln \ep|-b)\,,
\]
and evaluate $\ue(x,t^{\,\ep}(b))$ at a point $x \in \om _0
^{(1)}$ such that $dist(x,\Gamma _0)=\mathcal C\ep$. Define $C:=
\|u_0\|_{{C^2}(\overline{\Omega})}$. Since $u_0=a$ on $\Gamma_0$,
we have that
\begin{equation}\label{m-M}
u_0(x) \leq a +C \mathcal C \ep\,,
\end{equation}
which implies together with Lemma \ref{lemme-opt} that
\[
\begin{array}{lll}
w^+ _\ep (x, t^{\,\ep}(b))& =Y\Big ( \mu ^{-1} (|\ln \ep|-b),
u_0(x)+\ep \CG e^{-b}-\ep ^2 \CG \Big ) \vsp \\
& \leq a+C_2e^{|\ln \ep|-b} \big(u_0(x)+\ep \CG e^{-b}-\ep^2
\CG -a ) \vsp \\
& \leq a+C_2 \ep^{-1} e^{-b} (C\mathcal C \ep+\ep \CG e^{-b}) \vsp \\
& = a+C_2 e^{-b} (C\mathcal C + \CG e^{-b})\,.
\end{array}
\]
Now choose $b$ large enough so that
\begin{equation}\label{1-eta}
a+C_2 e^{-b} (C\mathcal C + \CG e^{-b}) <1-\eta\,;
\end{equation}
the inequalities \eqref{1-eta} and \eqref{g-coincee1} then yield
\[
\ue(x,t^{\,\ep}(b))\leq w^+ _\ep(x, t^{\,\ep}(b))<1-\eta\,.
\]
Therefore \eqref{resultat} does not hold at $t=t^{\,\ep}(b)$, and
hence $t^{\,\ep}(b)<t^{\,\ep}_{min}$.  This completes the proof of
Proposition \ref{pr:optimal-time}. \qed

\end{document}